\documentclass[pamm,a4paper,fleqn]{w-art}
\usepackage{times,cite,w-thm}
\usepackage[T1]{fontenc}
\usepackage[utf8]{inputenc}
\theoremstyle{plain}

\theoremstyle{definition}

\usepackage{graphicx}

\usepackage{pgfplotstable}
\usepackage{placeins}

\usepackage{amsmath,amssymb,amsthm}
\usepackage{xparse}

\newtheorem{algorithm}{Algorithm}

\let\epsilon\varepsilon
\let\phi\varphi
\let\rho\varrho

\usepackage{xspace}

\newcommand{\etal}[0]{{et~al.\@}\xspace}

\renewcommand{\O}{\mathcal{O}}

\newcommand{\N}{\mathbb{N}}
\providecommand\R{}% does nothing, if already defined
\renewcommand{\R}{\mathbb{R}}
\NewDocumentCommand{\RK}{o m O{\the\numexpr#2-1\relax} m O{} O{} o}{%
  \IfValueTF{#1}{#1}{RK}%
  #2(#3)#4%
  \ifblank{#6}{}{\textsubscript{F}}%
  \ifblank{#5}{}{[#5]}%
  \IfValueT{#7}{#7}%
}

\newcommand{\dt}{\Delta t}

\renewcommand{\vec}[1]{\pmb{#1}}

\NewDocumentCommand{\opD}{m+g}{%
  \IfNoValueTF{#2}
    {D_{#1}}
    {D_{#1,#2}}%
}
\NewDocumentCommand{\opDsplit}{m+g}{%
  \IfNoValueTF{#2}
    {\widetilde{D}_{#1}}
    {\widetilde{D}_{#1,#2}}%
}
\NewDocumentCommand{\opM}{g}{%
  \IfNoValueTF{#1}
    {M}
    {M_{#1}}%
}
\NewDocumentCommand{\opQ}{g}{%
  \IfNoValueTF{#1}
    {Q}
    {Q_{#1}}%
}
\NewDocumentCommand{\opI}{g}{%
  \IfNoValueTF{#1}
    {I}
    {I_{#1}}%
}
\NewDocumentCommand{\opV}{g}{%
  \IfNoValueTF{#1}
    {V}
    {V_{#1}}%
}
\NewDocumentCommand{\opB}{g}{%
  \IfNoValueTF{#1}
    {B}
    {B_{#1}}%
}
\NewDocumentCommand{\opR}{g}{%
  \IfNoValueTF{#1}
    {R}
    {R_{#1}}%
}
\NewDocumentCommand{\opN}{m+g}{%
  \IfNoValueTF{#2}
    {N_{#1}}
    {N_{#1,#2}}%
}
\NewDocumentCommand{\fnum}{g}{%
  \IfNoValueTF{#1}
    {f^{\mathrm{num}}}
    {f^{\mathrm{num,#1}}}%
}
\NewDocumentCommand{\vecfnum}{g}{%
  \IfNoValueTF{#1}
    {\vec{f}^{\mathrm{num}}}
    {\vec{f}^{\mathrm{num,#1}}}%
}
\NewDocumentCommand{\vecfcorr}{g}{%
  \IfNoValueTF{#1}
    {\vec{f}^{\mathrm{corr}}}
    {\vec{f}^{\mathrm{corr,#1}}}%
}
\NewDocumentCommand{\fvol}{g}{%
  \IfNoValueTF{#1}
    {f^{\smash{\mathrm{vol}}}}
    {f^{\smash{\mathrm{vol,#1}}}}%
}

\newcommand{\Te}{T_{end}}

\newcommand{\km}{k_{\max}}

\newcommand{\II}{\mathcal I}
\newcommand{\algs}[4]{#1^{#2,[#3],{#4}}}

\definecolor{col1}{RGB}{         0  135.4688  255.0000}
\definecolor{col2}{RGB}{         0    9.9609  255.0000}
\definecolor{col3}{RGB}{  115.5469         0  255.0000}
\definecolor{col4}{RGB}{  241.0547         0  255.0000}
\definecolor{col5}{RGB}{  255.0000         0  143.4375}
\definecolor{col6}{RGB}{  255.0000         0   17.9297}

\pgfplotscreateplotcyclelist{hierarchy6}{%
	{col1,thick,mark=o},
	{col2,thick,mark=pentagon},
	{col3,thick,mark=square},
	{col4,thick,mark=triangle},
	{col5,thick,mark=x},
   {col6,thick,mark=otimes}}

\begin{document}
%% \def\leftmark{Session title}
%%
%%    The information for the title page will be placed between
%%    \begin{document} and \maketitle. The order of most entries
%%    is determined by the class file and can not be changed by
%%    rearranging them. The maketitle command follows after the
%%    abstract.
%%
%%    Most of the following commands will be completed by the publisher.
%%
%%    \renewcommand{\copyrightyear}{2016}
%%    \DOIsuffix{pamm.20161zzzz}
%%    \Volume{16}
%%    \Year{2016}
%%    \pagespan{1}{}
%%
%%    The short title is optional:

\TitleLanguage[EN]
\title[Functional-preserving multiderivative schemes]{Functional-preserving predictor-corrector multiderivative schemes}

%% Please do not enter footnotes or \inst{}-notes into the optional
%% argument of the author command.

%% Please delete not needed author entries.
%% Information for the first author.
\author{\firstname{Hendrik}  \lastname{Ranocha}\inst{1,}%
  \footnote{\ElectronicMail{mail@ranocha.de}}}

\address[\inst{1}]{\CountryCode[DE]Applied Mathematics, University of Hamburg, Bundesstr. 55, 20146 Hamburg, Germany}
%%
%%    Information for the second author
\author{\firstname{Jochen} \lastname{Sch\"utz}\inst{2,}%
  \footnote{Corresponding author: \ElectronicMail{jochen.schuetz@uhasselt.be}}}

%%    Information for the third author
\author{\firstname{Eleni} \lastname{Theodosiou}\inst{2,}%
  \footnote{\ElectronicMail{eleni.theodosiou@uhasselt.be}}}

\address[\inst{2}]{\CountryCode[BE]Faculty of Sciences \& Data Science Institute, Hasselt University, Agoralaan Gebouw D, 3590 Diepenbeek, Belgium}
%%

%%    Abstract is required.
\AbstractLanguage[EN]
\begin{abstract}
  In this work, we develop a class of high-order multiderivative time integration methods that is able to preserve certain functionals discretely. Important ingredients are the recently developed Hermite-Birkhoff-Predictor-Corrector methods and the technique of relaxation for numerical methods of ODEs. We explain the algorithm in detail and show numerical results for two- and three-derivative methods, comparing relaxed and unrelaxed methods. The numerical results demonstrate that, at the slight cost of the relaxation, an improved scheme is obtained.
\end{abstract}
\maketitle

%%%%%%%%%%%%%%%%%%%%%%%%%%%%%%

\section{Introduction}\label{sec:intro}

The efficient and accurate numerical solution of time-dependent differential equations is ubiquitous in the computational sciences; examples of practical interest stem from meteorology, aerospace engineering, porous media flow and many more.
There are several challenges associated to high-order temporal integration, such as efficiency and stability, which are obviously intertwined.
In classical numerical schemes, high-order has been reached through an increase in either stages or steps, or both, see, e.g., \cite{butcher2006general}. By now, Runge-Kutta schemes and linear multistep schemes are a de-facto standard in, e.g., the computational fluid dynamics community (CFD), see \cite{Hartmann2021} for an overview on the use of implicit methods in CFD. Although also a rather classical approach, see \cite{HaWa73}, the multiderivative paradigm has only been rediscovered rather recently; for some examples see  \cite{TC10,Seal13,SSJ2017,2022_Gottlieb_EtAl,MoradiAbdi2022,Zeifang2021} and the references therein.

To illustrate the approach, let us assume that the underlying differential equation is given by
\begin{alignat}{2}
 \label{eq:ode}
 w'(t) &= \Phi(w(t)), \quad t \in [0, \Te], \qquad  w(0) &\ =\ & w_0,
\end{alignat}
for some unknown function $w : [0, \Te] \rightarrow \R^{\dim}$ and a given smooth function $\Phi: \R^{\dim} \rightarrow \R^{\dim}$. Obviously, the second derivative of $w$ can be computed from $\Phi$ and its Jacobian through
\begin{align}
 \label{eq:dotf}
 w''(t) = \Phi'(w(t)) \Phi(w(t)) =: \dot \Phi(w(t)).
\end{align}
Obviously, also the third temporal derivative $\ddot \Phi(w)$ and higher derivatives of $w$ can be computed.
Multiderivative time integrators explicitly take the quantities $\Phi, \dot \Phi, \ldots$ into account, which results for, e.g., a given number of stages, in a higher order than in a classical approach.
In this work, we consider a peculiar predictor-corrector form of the implicit multiderivative method, inspired by spectrally deferred correction methods \cite{OngSpiteri}. This HBPC (Hermite-Birkhoff-Predictor-Corrector) method was initially developed and motivated as an IMEX scheme in \cite{SealSchuetz19} and then subsequently extended to higher orders in \cite{SchuetzSealZeifang21,ZSS22}.
HBPC has shown favorable behavior for the solution of compressible flow equations \cite{Zeifang2021,ZMS2022,Arjun2023}.

While linear stability of HBPC has been tackled in \cite{ZSS22}, the behaviour of the method for large values of $\Te$ has not been investigated yet. As for most schemes, it is to be expected that the numerical error grows tremendously with growing $\Te$. In this work, we consider the case of a functional $\eta: \R^{\dim} \rightarrow \R$ that is preserved under the solution, i.e.,
\begin{align}
 \label{eq:conservedfunctional}
 \frac{\mathrm{d}}{\mathrm d t} \eta(w(t)) \equiv 0.
\end{align}
For Hamiltonian problems, $\eta$ could simply be the Hamiltonian function; for smooth flow problems, it could be entropy and so on.
First, we show how the classical HBPC method behaves in terms of $\eta$ and in terms of the numerical error growth over time. Second, we extend the HBPC method with a relaxation procedure, originally developed in \cite{ketcheson2019relaxation,ranocha2020relaxation,ranocha2020general}; based on an older idea from \cite{sanzserna1982explicit}.
This relaxation procedure, outlined below, enforces the preservation of $\eta$ through an additional projection step.
This projection step necessitates the solution of a scalar equation, typically through Newton's method or more efficient variants of the bisection method. While for explicit low-order schemes, this might constitute a significant overhead \cite{ranocha2020relaxationHamiltonian}, it is negligible in our setting of implicit schemes. We show that with this very simple addendum to the algorithm, both error growth in time is reduced and the functional $\eta$ is preserved for several testcases.

%The paper is outlined as follows: \todo{todo!}

\section{Numerical tools}\label{sec:numerics}
In this chapter, we describe the underlying time integration algorithm as well as its combination with relaxation.
In the following, $t^n$ refers to the time instance $t^n := n \dt$ with some fixed (only for the ease of presentation) timestep $\dt > 0$.
%Please note that a fixed $\dt$ is not a necessity for the algorithms to work.

\subsection{Hermite-Birkhoff predictor-corrector time integration}
The algorithm to be explained in the following is of the predictor-corrector type, iterating towards a background, fully implicit multiderivative Runge-Kutta scheme using $m\in\N$ temporal derivatives of $w$.
For the ease of presentation, we first define this background scheme. Please note that this scheme is not actually used in our computations, only through the use of the corresponding quadrature rule.
The scheme is of the classical multiderivative Runge-Kutta type, with $s$ stages $w^{n, l}$, $1 \leq l \leq s$, and update $w^{n+1}_{RK}$ defined by:
\begin{equation}
\label{eq:rk}
\begin{aligned}
 w^{n,l} := w^n + \sum_{d=1}^m \dt^d \sum_{j=1}^s B_{lj}^{(d)} \frac{\mathrm d^{d-1}}{\mathrm d t^{d-1}} \Phi(w^{n,j}), \quad
 w^{n+1}_{RK} := w^n + \sum_{d=1}^m \dt^d \sum_{j=1}^s b_{j}^{(d)} \frac{\mathrm d^{d-1}}{\mathrm d t^{d-1}} \Phi(w^{n,j}).
\end{aligned}
\end{equation}
The matrices $B^{(d)}$, $1 \leq d \leq m$, form the Butcher tableaux. It is assumed that the $l-$th stage value of time is $t^n + c_l \dt$, for values $c_l \equiv \sum_{j=1}^s B_{lj}^{(1)}$. The coefficients for the Runge-Kutta update are denoted by $b^{(d)}_l, 1 \leq l \leq s$.
We assume that the Runge-Kutta scheme associated with this Butcher tableau is of order $q$.
Please note that we have defined
\begin{align*}
\frac{\mathrm d^0}{\mathrm d t^0} \Phi(w^{n,j}) := \Phi(w^{n,j}), \quad \frac{\mathrm d^1}{\mathrm d t^1} \Phi(w^{n,j}):= \dot\Phi(w^{n,j}),
\quad \frac{\mathrm d^2}{\mathrm d t^2} \Phi(w^{n,j}) := \ddot \Phi(w^{n,j}),
\end{align*}
and so on.
In this work, we rely on three Runge-Kutta schemes, two with two-derivatives, see \cite[Eq.~(2) and Eq. (3), respectively, for the Butcher tableaux]{SchuetzSealZeifang21}, and a two-point three-derivative scheme with Butcher tableau
\begin{align*}
 c = \begin{pmatrix} 0 \\ 1 \end{pmatrix}, \quad B^{(1)} = \begin{pmatrix} 0 & 0 \\ \frac 1 2 & \frac 1 2 \end{pmatrix}, \quad
												 B^{(2)} = \begin{pmatrix} 0 & 0 \\ \frac 1 {10} &- \frac 1 {10} \end{pmatrix}, \quad
												 B^{(3)} = \begin{pmatrix} 0 & 0 \\ \frac 1 {120} & \frac 1 {120} \end{pmatrix}.
\end{align*}
The final HBPC scheme to be presented here relies on a predictor ($k=0$) and correction steps ($1 \leq k \leq \km$) for the quantities $w^{n,l}$. For short, the notation here is $w^{n, [k], l}$. The predictor is a straightforward implicit Taylor scheme making use of $m$ temporal derivatives of $w$, the corrector is very similar in structure plus it additionally relies on the quadrature formula $\II_l$ defined through the Runge-Kutta scheme \eqref{eq:rk} by
\begin{align*}
 	\II_l:=&
 	\sum_{d=1}^m \dt^d \sum_{j=1}^s B^{(d)}_{lj}  \algs{\frac{\mathrm{d}^{d-1}}{\mathrm{d}t^{d-1}}\Phi}{n}{k}{j}.
 \end{align*}
Note the shorthand notation $\algs\Phi n k j := \Phi(\algs w n k j)$.
Finally, we obtain
\begin{algorithm}[HBPC($m$, $q$, $\km$) \cite{SealSchuetz19,SchuetzSealZeifang21}]\label{alg:mdode_parallel}
    The algorithm consists of the following three steps:
	\begin{enumerate}
 		\item \textbf{Predict.} Solve the following expression for $\algs w n 0 l$ and $1 \leq l \leq s$:
 		\begin{align}
 		\label{eq:predictor_Taylor}
 		\begin{split}
 		\algs w {n} 0 l := w^{n} + \sum_{d = 1}^m \frac{(-1)^{d-1}(c_l \dt)^d}{d!} \algs {\frac{\mathrm{d}^{d-1}}{\mathrm{d}t^{d-1}}\Phi} n 0 l .
 		\end{split}
 		\end{align}
 		Subsequently:
 		\item \textbf{Correct.} Solve the following for $\algs w {n} {k+1} l$, for each $1 \leq l \leq s$ and each $0 \leq k < k_{\max}$:
 		\begin{align}
 		\begin{split}
 		\label{eq:correct_parallel}
 		\algs w {n} {k+1} l &:= w^n + \sum_{d=1}^{m} \frac{(-1)^{d-1}\dt^d}{d!} \left( \algs {\frac{\mathrm{d}^{d-1}}{\mathrm{d}t^{d-1}}\Phi} {n} {k+1} l - \algs {\frac{\mathrm{d}^{d-1}}{\mathrm{d}t^{d-1}}\Phi} n {k} l\right) +
 		\II_l
 		\end{split}
 		\end{align}
 		\item \textbf{Update.} Set
 		\begin{align}
 		 \label{eq:update}
 		 w^{n+1} := w^n + \sum_{d=1}^{m} \frac{(-1)^{d-1}\dt^d}{d!} \frac{\mathrm{d}^{d-1}}{\mathrm{d}t^{d-1}} \left( \algs {\Phi} {n} {\km} l - \algs {\Phi} n {\km-1} l\right) + \sum_{d=1}^m \dt^d \sum_{j=1}^s b^{(d)}_j \algs{\frac{\mathrm{d}^{d-1}}{\mathrm{d}t^{d-1}}\Phi}{n}{\km-1}{j}.
 		\end{align}
 	\end{enumerate}
\end{algorithm}

\begin{remark}
 Please note that whenever the background Runge-Kutta scheme is globally stiffly accurate, i.e., there holds
 \begin{align*}
  b_j^{(d)} = B^{(d)}_{sj}, \quad 1 \leq j \leq s, 1 \leq d \leq m,
 \end{align*}
 then the update step reduces to $w^{n+1} := \algs w n \km s$.
 It is hence a slight generalization of \cite{SealSchuetz19,SchuetzSealZeifang21}, where only schemes with $c_s = 1$ are treated.
 In any case, the update step is explicit.
\end{remark}

\begin{remark}\label{rem:order}
 The order of convergence $p$ of this scheme is the minimum of $k_{\max}+m$ and the order $q$ of the underlying Runge-Kutta scheme; hence, $p := \min\{k_{max} + m, q\}$.
\end{remark}

\subsection{Relaxation procedure}

The idea of a relaxation procedure as introduced in \cite{ketcheson2019relaxation,ranocha2020relaxation,ranocha2020general} is to consider a \emph{scalar} parameter $\gamma \in \R$ and form a linear combination of $w^n$ and $w^{n+1}$ to obtain the quantity $w^{n+1}_\gamma = w^n + \gamma (w^{n+1}-w^n)$. The relaxation parameter $\gamma$ gives the flexibility to enforce the preservation of the functional $\eta$, just as for the continuous case, see Eq.~\eqref{eq:conservedfunctional}, through the equation (in $\gamma$)
\begin{align}
 \label{eq:relaxation}
 \eta\left(w^n + \gamma (w^{n+1}-w^n)\right) = \eta(w^n).
\end{align}
After having found a suitable $\gamma$ -- typically through a scalar Newton algorithm --, the relaxed update $w^{n+1}_\gamma = w^n + \gamma (w^{n+1}-w^n)$ is considered the new update step at time level $t^n + \gamma \dt$. Computation with Alg.~\ref{alg:mdode_parallel} then continues from this adapted point in time and the corresponding linear combination of $w^n$ and $w^{n+1}$ as usual. Note that $\dt$ is a constant throughout the computation (although adaptive timesteps are a possibility as well), however, the resulting time instances are not necessarily spaced equidistantly.

Obviously, $\gamma = 0$ is a (meaningless) solution to \eqref{eq:relaxation}. It has been shown in \cite{ranocha2020general} that under rather mild conditions on the timestep $\dt$ and the functional $\eta$, there is also a unique solution $\gamma$ which is close to one, in fact, it is $\O(\dt^{p+1}$) away from one. Here, $p$ denotes the order of the method. With this solution, the relaxation approach keeps at least the order of the baseline methods. The relaxation approach is not restricted to invariants and has also been extended to general functionals $\eta$ in \cite{ketcheson2019relaxation,ranocha2020relaxation,ranocha2020general}, resulting for example in efficient, fully-discrete, and locally entropy-stable numerical methods for computational fluid dynamics \cite{ranocha2020fully}
and nonlinear dispersive wave equations \cite{ranocha2021broad,mitsotakis2021conservative,ranocha2021rate}.

\section{Numerical experiments}\label{sec:numres}
In this section, we present numerical findings of the HBPC method for a couple of test problems. As we are dealing with implicit time integration, both linear and nonlinear solvers are important ingredients. In all the numerical results to follow, we use a damped Newton procedure for the nonlinear equations, together with the standard backslash operator in Matlab to solve the linear systems. The Newton tolerance is always set to a very fine tolerance $10^{-14}$, and a maximum of 1000 iterations is allowed. Obviously, we did not go for the most efficient solution here. For considerations regarding Newton efficiency, we refer the reader to \cite{ZMS2022}. In all the numerical results, 'error' is defined as the Euclidean error of the discrete solution at the final time $\Te$.

\subsection{Nonlinear oscillator}

As a first numerical example, we consider the nonlinear oscillator of \cite{ranocha2021strong,ranocha2020energy}, given by
\begin{align*}
 \Phi(w) := \frac{1}{\|w\|_2^2} \begin{pmatrix} -w_2 \\ w_1 \end{pmatrix}, \qquad w(0) = \begin{pmatrix} 1 \\ 0 \end{pmatrix}.
\end{align*}
The standard squared Euclidean norm is a conservative functional for this problem, i.e., $\eta(w) := \|w\|_2^2$ is a constant along the solution for all times $t \in \R^+$.

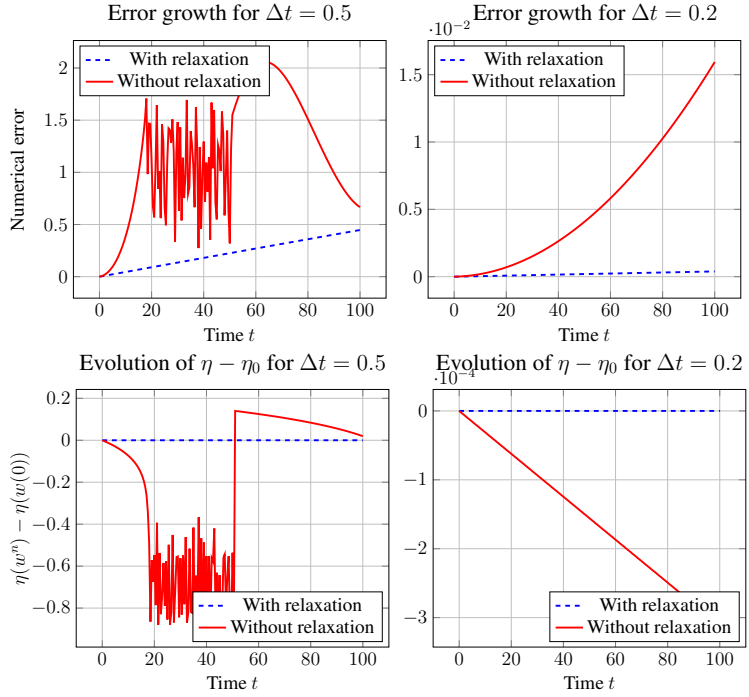
\begin{figure} [ht]
        \centering
        \sidecaption
        \begin{minipage}{0.6\textwidth}
        \begin{tikzpicture}[scale=0.6]
            \begin{axis}[
              xlabel={Time $t$},
              ylabel={Numerical error},
              label style={font=\large},
              grid=major,
              legend style={at={(0.02,0.98)},anchor= north west,font=\large},
              title={Error growth for $\dt = 0.5$},
              title style={font=\Large},
              xticklabel style = {font=\large},
              yticklabel style = {font=\large},
            ]
            \addplot +[very thick,mark=none,dashed] table [x expr={\thisrowno{0}}, y expr={\thisrowno{1}}]
                {numerical_results/Oscillator/EntropyConservative/ErrorGrowth_HBPC_serial_HB-I2DRK6-3s_relaxation=1_kmax=4_Predictor=Taylor_oscillator_Tend=100_nonlin=1_damped_newton_backslash_Nt=200_norm=l2.csv};
            \addplot +[very thick,mark=none] table [x expr={\thisrowno{0}}, y expr={\thisrowno{1}}]
                {numerical_results/Oscillator/EntropyConservative/ErrorGrowth_HBPC_serial_HB-I2DRK6-3s_relaxation=0_kmax=4_Predictor=Taylor_oscillator_Tend=100_nonlin=1_damped_newton_backslash_Nt=200_norm=l2.csv};
              \legend{{With relaxation}, {Without relaxation}}
            \end{axis}
            \end{tikzpicture}
        \begin{tikzpicture}[scale=0.6]
            \begin{axis}[
              xlabel={Time $t$},
              label style={font=\large},
              grid=major,
              legend style={at={(0.02,0.98)},anchor= north west,font=\large},
              title={Error growth for $\dt = 0.2$},
              title style={font=\Large, xshift=1.5ex},
              xticklabel style = {font=\large},
              yticklabel style = {font=\large},
            ]
            \addplot +[very thick,mark=none,dashed] table [x expr={\thisrowno{0}}, y expr={\thisrowno{1}}]
                {numerical_results/Oscillator/EntropyConservative/ErrorGrowth_HBPC_serial_HB-I2DRK6-3s_relaxation=1_kmax=4_Predictor=Taylor_oscillator_Tend=100_nonlin=1_damped_newton_backslash_Nt=500_norm=l2.csv};
            \addplot +[very thick,mark=none] table [x expr={\thisrowno{0}}, y expr={\thisrowno{1}}]
                {numerical_results/Oscillator/EntropyConservative/ErrorGrowth_HBPC_serial_HB-I2DRK6-3s_relaxation=0_kmax=4_Predictor=Taylor_oscillator_Tend=100_nonlin=1_damped_newton_backslash_Nt=500_norm=l2.csv};
                \legend{{With relaxation}, {Without relaxation}}
            \end{axis}
          \end{tikzpicture} \\
          \begin{tikzpicture}[scale=0.6]
            \begin{axis}[
              xlabel={Time $t$},
              ylabel={$\eta(w^n)-\eta(w(0))$},
              label style={font=\large},
              grid=major,
              legend style={at={(0.98,0.02)},anchor= south east,font=\large},
              title={Evolution of $\eta-\eta_0$ for $\dt = 0.5$},
              title style={font=\Large},
              xticklabel style = {font=\large},
              yticklabel style = {font=\large},
            ]
            \addplot +[very thick,mark=none,dashed] table [x expr={\thisrowno{0}}, y expr={\thisrowno{2}-1}]
                {numerical_results/Oscillator/EntropyConservative/ErrorGrowth_HBPC_serial_HB-I2DRK6-3s_relaxation=1_kmax=4_Predictor=Taylor_oscillator_Tend=100_nonlin=1_damped_newton_backslash_Nt=200_norm=l2.csv};
            \addplot +[very thick,mark=none] table [x expr={\thisrowno{0}}, y expr={\thisrowno{2}-1}]
                {numerical_results/Oscillator/EntropyConservative/ErrorGrowth_HBPC_serial_HB-I2DRK6-3s_relaxation=0_kmax=4_Predictor=Taylor_oscillator_Tend=100_nonlin=1_damped_newton_backslash_Nt=200_norm=l2.csv};
              \legend{{With relaxation}, {Without relaxation}}
            \end{axis}
          \end{tikzpicture}
        \begin{tikzpicture}[scale=0.6]
            \begin{axis}[
              xlabel={Time $t$},%ylabel={$\eta(w^n)$},
              label style={font=\large},
              grid=major,
              legend style={at={(0.98,0.02)},anchor= south east,font=\large},
              title={Evolution of $\eta-\eta_0$ for $\dt = 0.2$},
              title style={font=\Large},
              xticklabel style = {font=\large},
              yticklabel style={/pgf/number format/fixed,
                  /pgf/number format/precision=6,font=\large}
            ]%,ymax=1e1,ymin=1e-13]
            \addplot +[very thick,mark=none,dashed] table [x expr={\thisrowno{0}}, y expr={\thisrowno{2}-1}]
                {numerical_results/Oscillator/EntropyConservative/ErrorGrowth_HBPC_serial_HB-I2DRK6-3s_relaxation=1_kmax=4_Predictor=Taylor_oscillator_Tend=100_nonlin=1_damped_newton_backslash_Nt=500_norm=l2.csv};
            \addplot +[very thick,mark=none] table [x expr={\thisrowno{0}}, y expr={\thisrowno{2}-1}]
                {numerical_results/Oscillator/EntropyConservative/ErrorGrowth_HBPC_serial_HB-I2DRK6-3s_relaxation=0_kmax=4_Predictor=Taylor_oscillator_Tend=100_nonlin=1_damped_newton_backslash_Nt=500_norm=l2.csv};
                \legend{{With relaxation}, {Without relaxation}}
            \end{axis}
          \end{tikzpicture}
          \end{minipage}
        \caption{The (serial) HBPC(2, 6, 4) predictor-corrector scheme of \cite{ZSS22}, see also \cite{SchuetzSealZeifang21}, applied to the entropy-conserving nonlinear oscillator with $\Te = 100$. Left are numerical results for $\dt = 0.5$, right are results for $\dt = 0.2$. Top: numerical error as a function of time; bottom: the deviation in the functional $\eta$ evaluated for the discrete solution. It is clearly visible that in all cases, the relaxed method behaves significantly better than its unrelaxed counterpart.
        }\label{fig:ErrorGrowthOscillator}
    \end{figure}

\paragraph{Error growth}
In a first step, we consider the error growth for the HBPC scheme in dependency of time with and without relaxation. As final time, the rather large $\Te = 100$ is chosen in combination with the large timesteps $\dt = 0.5$ and $\dt = 0.2$, respectively.
As a time integrator, the HBPC(2,6,4) method is used, i.e., order six is to be expected. Please note that the behavior of this method is representative.
Time against error can be seen in the top of Fig.~\ref{fig:ErrorGrowthOscillator} for the algorithm with and without relaxation. It can be clearly seen that the numerical error for the relaxed HBPC method behaves linearly in both cases. At least for smaller $t$, the error of the unrelaxed method behaves quadratically. For $\dt = 0.5$, it starts to oscillate at some point. This is also reflected in the fact that Newton's algorithm did not converge for the unrelaxed method and $\dt = 0.5$. In this sense, the relaxation improved the algorithm tremendously, even if one is not interested in an accurate representation of $\eta$. The bottom of Fig.~\ref{fig:ErrorGrowthOscillator} shows the evolution of $\eta-\eta_0$ ($\eta_0 := \eta(w(0))$) for the two values of $\dt$. As expected, the relaxed version preserves $\eta$, even if the error level, at least for $\dt = 0.5$, is also rather high for the relaxed method.
All these results are very much in line with the results from literature as presented in \cite{cano1997error,duran1998numerical,calvo2011error}.

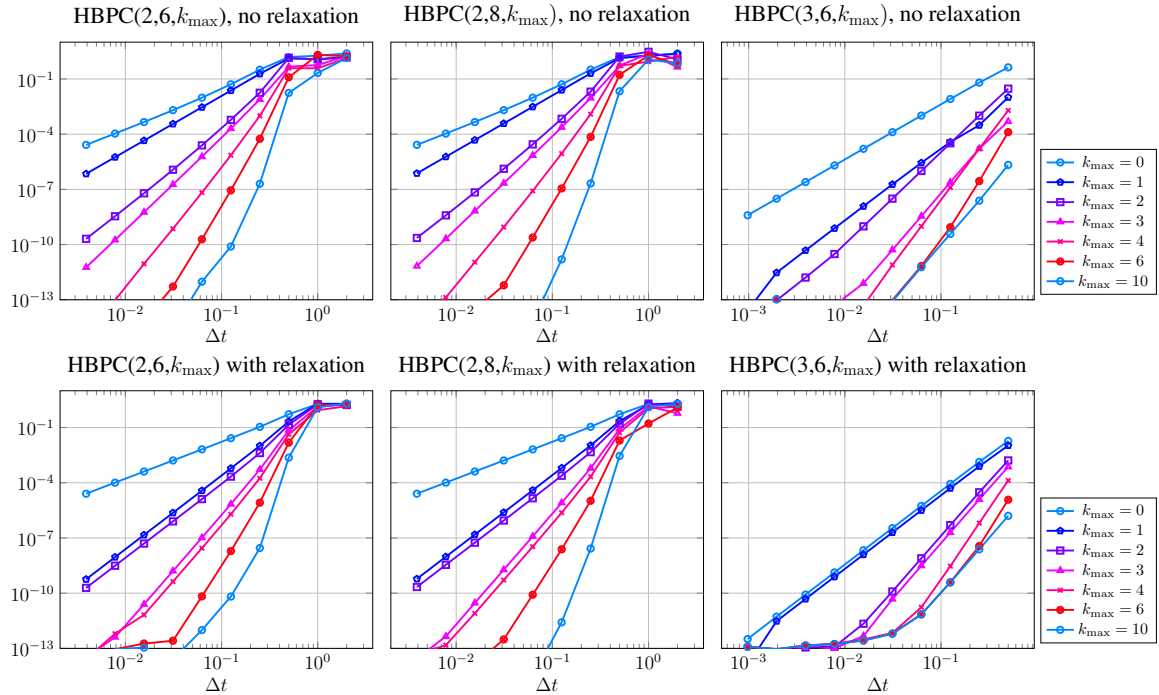
\begin{figure} [ht]
        %\centering
        %\sidecaption
        %\begin{minipage}{0.5\textwidth}
        \begin{tikzpicture}[scale=0.6]
            \begin{loglogaxis}[
            cycle list name=hierarchy6,
            xlabel={$\Delta t$},
            ylabel={},
            grid=major,
            title={HBPC(2,6,$k_{\max}$), no relaxation},
            label style={font=\large},
            title style={font=\Large},
            ymax=1e1,ymin=1e-13,
            xticklabel style = {font=\large},
            yticklabel style = {font=\large},
          ]
            \addplot +[very thick] table [x expr={\thisrowno{0}}, y expr={\thisrowno{1}},
                            skip first n=1, col sep=comma]
                {numerical_results/Oscillator/EntropyConservative/Results_HBPC_serial_HB-I2DRK6-3s_relaxation=0_kmax=0_Predictor=Taylor_oscillator_Tend=10_nonlin=1_damped_newton_backslash_norm=l2.csv};
            \addplot +[very thick] table [x expr={\thisrowno{0}}, y expr={\thisrowno{1}},
                            skip first n=1, col sep=comma]
                {numerical_results/Oscillator/EntropyConservative/Results_HBPC_serial_HB-I2DRK6-3s_relaxation=0_kmax=1_Predictor=Taylor_oscillator_Tend=10_nonlin=1_damped_newton_backslash_norm=l2.csv};
                \addplot +[very thick] table [x expr={\thisrowno{0}}, y expr={\thisrowno{1}},
                            skip first n=1, col sep=comma]
                {numerical_results/Oscillator/EntropyConservative/Results_HBPC_serial_HB-I2DRK6-3s_relaxation=0_kmax=2_Predictor=Taylor_oscillator_Tend=10_nonlin=1_damped_newton_backslash_norm=l2.csv};
                \addplot +[very thick] table [x expr={\thisrowno{0}}, y expr={\thisrowno{1}},
                            skip first n=1, col sep=comma]
                {numerical_results/Oscillator/EntropyConservative/Results_HBPC_serial_HB-I2DRK6-3s_relaxation=0_kmax=3_Predictor=Taylor_oscillator_Tend=10_nonlin=1_damped_newton_backslash_norm=l2.csv};
                \addplot +[very thick] table [x expr={\thisrowno{0}}, y expr={\thisrowno{1}},
                            skip first n=1, col sep=comma]
                {numerical_results/Oscillator/EntropyConservative/Results_HBPC_serial_HB-I2DRK6-3s_relaxation=0_kmax=4_Predictor=Taylor_oscillator_Tend=10_nonlin=1_damped_newton_backslash_norm=l2.csv};
                \addplot +[very thick] table [x expr={\thisrowno{0}}, y expr={\thisrowno{1}},
                            skip first n=1, col sep=comma]
                {numerical_results/Oscillator/EntropyConservative/Results_HBPC_serial_HB-I2DRK6-3s_relaxation=0_kmax=6_Predictor=Taylor_oscillator_Tend=10_nonlin=1_damped_newton_backslash_norm=l2.csv};\addplot +[very thick] table [x expr={\thisrowno{0}}, y expr={\thisrowno{1}},
                            skip first n=1, col sep=comma]
                {numerical_results/Oscillator/EntropyConservative/Results_HBPC_serial_HB-I2DRK6-3s_relaxation=0_kmax=10_Predictor=Taylor_oscillator_Tend=10_nonlin=1_damped_newton_backslash_norm=l2.csv};
            %\legend{$k_{\max}=0$,$k_{\max}=1$,$k_{\max}=2$,$k_{\max}=3$,$k_{\max}=4$,$k_{\max}=10$}
            \end{loglogaxis}
            \end{tikzpicture}
            \begin{tikzpicture}[scale=0.6]
            \begin{loglogaxis}[
              cycle list name=hierarchy6,
              xlabel={$\Delta t$},
              ylabel={},
              grid=major,
              title={HBPC(2,8,$k_{\max})$, no relaxation},
              label style={font=\large},
              title style={font=\Large},
              ymax=1e1,ymin=1e-13,
              yticklabels=\empty,
              xticklabel style = {font=\large},
              yticklabel style = {font=\large},
            ]
            \addplot +[very thick] table [x expr={\thisrowno{0}}, y expr={\thisrowno{1}},
                            skip first n=1, col sep=comma]
                {numerical_results/Oscillator/EntropyConservative/Results_HBPC_serial_HB-I2DRK8-4s_relaxation=0_kmax=0_Predictor=Taylor_oscillator_Tend=10_nonlin=1_damped_newton_backslash_norm=l2.csv};
            \addplot +[very thick] table [x expr={\thisrowno{0}}, y expr={\thisrowno{1}},
                            skip first n=1, col sep=comma]
                {numerical_results/Oscillator/EntropyConservative/Results_HBPC_serial_HB-I2DRK8-4s_relaxation=0_kmax=1_Predictor=Taylor_oscillator_Tend=10_nonlin=1_damped_newton_backslash_norm=l2.csv};
                \addplot +[very thick] table [x expr={\thisrowno{0}}, y expr={\thisrowno{1}},
                            skip first n=1, col sep=comma]
                {numerical_results/Oscillator/EntropyConservative/Results_HBPC_serial_HB-I2DRK8-4s_relaxation=0_kmax=2_Predictor=Taylor_oscillator_Tend=10_nonlin=1_damped_newton_backslash_norm=l2.csv};
                \addplot +[very thick] table [x expr={\thisrowno{0}}, y expr={\thisrowno{1}},
                            skip first n=1, col sep=comma]
                {numerical_results/Oscillator/EntropyConservative/Results_HBPC_serial_HB-I2DRK8-4s_relaxation=0_kmax=3_Predictor=Taylor_oscillator_Tend=10_nonlin=1_damped_newton_backslash_norm=l2.csv};
                \addplot +[very thick] table [x expr={\thisrowno{0}}, y expr={\thisrowno{1}},
                            skip first n=1, col sep=comma]
                {numerical_results/Oscillator/EntropyConservative/Results_HBPC_serial_HB-I2DRK8-4s_relaxation=0_kmax=4_Predictor=Taylor_oscillator_Tend=10_nonlin=1_damped_newton_backslash_norm=l2.csv};
                \addplot +[very thick] table [x expr={\thisrowno{0}}, y expr={\thisrowno{1}},
                            skip first n=1, col sep=comma]
                {numerical_results/Oscillator/EntropyConservative/Results_HBPC_serial_HB-I2DRK8-4s_relaxation=0_kmax=6_Predictor=Taylor_oscillator_Tend=10_nonlin=1_damped_newton_backslash_norm=l2.csv};\addplot +[very thick] table [x expr={\thisrowno{0}}, y expr={\thisrowno{1}},
                            skip first n=1, col sep=comma]
                {numerical_results/Oscillator/EntropyConservative/Results_HBPC_serial_HB-I2DRK8-4s_relaxation=0_kmax=10_Predictor=Taylor_oscillator_Tend=10_nonlin=1_damped_newton_backslash_norm=l2.csv};
            %\legend{$k_{\max}=0$,$k_{\max}=1$,$k_{\max}=2$,$k_{\max}=3$,$k_{\max}=4$,$k_{\max}=10$}
            \end{loglogaxis}
            \end{tikzpicture}
                        \begin{tikzpicture}[scale=0.6]
            \begin{loglogaxis}[
              cycle list name=hierarchy6,
              xlabel={$\Delta t$},
              ylabel={},
              grid=major,
              yticklabels=\empty,
              legend style={at={(1.02,0.02)},anchor= south west,font=\normalsize},
              title={HBPC(3,6,$k_{\max})$, no relaxation},
              label style={font=\large},
              title style={font=\Large},
              ymax=1e1,ymin=1e-13,
              xticklabel style = {font=\large},
              yticklabel style = {font=\large}
            ]
            \addplot +[very thick] table [x expr={\thisrowno{0}}, y expr={\thisrowno{1}},
                            skip first n=1, col sep=comma]
                {numerical_results/Oscillator/EntropyConservative/Results_HBPC_serial_HB-I3DRK6-2s_relaxation=0_kmax=0_Predictor=Taylor_oscillator_Tend=10_nonlin=1_damped_newton_backslash_norm=l2.csv};
            \addplot +[very thick] table [x expr={\thisrowno{0}}, y expr={\thisrowno{1}},
                            skip first n=1, col sep=comma]
                {numerical_results/Oscillator/EntropyConservative/Results_HBPC_serial_HB-I3DRK6-2s_relaxation=0_kmax=1_Predictor=Taylor_oscillator_Tend=10_nonlin=1_damped_newton_backslash_norm=l2.csv};
                \addplot +[very thick] table [x expr={\thisrowno{0}}, y expr={\thisrowno{1}},
                            skip first n=1, col sep=comma]
                {numerical_results/Oscillator/EntropyConservative/Results_HBPC_serial_HB-I3DRK6-2s_relaxation=0_kmax=2_Predictor=Taylor_oscillator_Tend=10_nonlin=1_damped_newton_backslash_norm=l2.csv};
                \addplot +[very thick] table [x expr={\thisrowno{0}}, y expr={\thisrowno{1}},
                            skip first n=1, col sep=comma]
                {numerical_results/Oscillator/EntropyConservative/Results_HBPC_serial_HB-I3DRK6-2s_relaxation=0_kmax=3_Predictor=Taylor_oscillator_Tend=10_nonlin=1_damped_newton_backslash_norm=l2.csv};
                \addplot +[very thick] table [x expr={\thisrowno{0}}, y expr={\thisrowno{1}},
                            skip first n=1, col sep=comma]
                {numerical_results/Oscillator/EntropyConservative/Results_HBPC_serial_HB-I3DRK6-2s_relaxation=0_kmax=4_Predictor=Taylor_oscillator_Tend=10_nonlin=1_damped_newton_backslash_norm=l2.csv};
                \addplot +[very thick] table [x expr={\thisrowno{0}}, y expr={\thisrowno{1}},
                            skip first n=1, col sep=comma]
                {numerical_results/Oscillator/EntropyConservative/Results_HBPC_serial_HB-I3DRK6-2s_relaxation=0_kmax=6_Predictor=Taylor_oscillator_Tend=10_nonlin=1_damped_newton_backslash_norm=l2.csv};
                \addplot +[very thick] table [x expr={\thisrowno{0}}, y expr={\thisrowno{1}},
                            skip first n=1, col sep=comma]
                {numerical_results/Oscillator/EntropyConservative/Results_HBPC_serial_HB-I3DRK6-2s_relaxation=0_kmax=10_Predictor=Taylor_oscillator_Tend=10_nonlin=1_damped_newton_backslash_norm=l2.csv};
            \legend{$k_{\max}=0$,$k_{\max}=1$,$k_{\max}=2$,$k_{\max}=3$,$k_{\max}=4$,$k_{\max}=6$,$k_{\max}=10$}
            \end{loglogaxis}
            \end{tikzpicture}\\
        \begin{tikzpicture}[scale=0.6]
            \begin{loglogaxis}[
            cycle list name=hierarchy6,
            xlabel={$\Delta t$},
            ylabel={},
            grid=major,
            legend style={at={(1.1,0.02)},anchor= south east,font=\large},
            title={HBPC(2,6,$k_{\max}$) with relaxation},
            label style={font=\large},
            title style={font=\Large},
            ymax=1e1,ymin=1e-13,
            xticklabel style = {font=\large},
            yticklabel style = {font=\large},
          ]
            \addplot +[very thick] table [x expr={\thisrowno{0}}, y expr={\thisrowno{1}},
                            skip first n=1, col sep=comma]
                {numerical_results/Oscillator/EntropyConservative/Results_HBPC_serial_HB-I2DRK6-3s_relaxation=1_kmax=0_Predictor=Taylor_oscillator_Tend=10_nonlin=1_damped_newton_backslash_norm=l2.csv};
            \addplot +[very thick] table [x expr={\thisrowno{0}}, y expr={\thisrowno{1}},
                            skip first n=1, col sep=comma]
                {numerical_results/Oscillator/EntropyConservative/Results_HBPC_serial_HB-I2DRK6-3s_relaxation=1_kmax=1_Predictor=Taylor_oscillator_Tend=10_nonlin=1_damped_newton_backslash_norm=l2.csv};
                \addplot +[very thick] table [x expr={\thisrowno{0}}, y expr={\thisrowno{1}},
                            skip first n=1, col sep=comma]
                {numerical_results/Oscillator/EntropyConservative/Results_HBPC_serial_HB-I2DRK6-3s_relaxation=1_kmax=2_Predictor=Taylor_oscillator_Tend=10_nonlin=1_damped_newton_backslash_norm=l2.csv};
                \addplot +[very thick] table [x expr={\thisrowno{0}}, y expr={\thisrowno{1}},
                            skip first n=1, col sep=comma]
                {numerical_results/Oscillator/EntropyConservative/Results_HBPC_serial_HB-I2DRK6-3s_relaxation=1_kmax=3_Predictor=Taylor_oscillator_Tend=10_nonlin=1_damped_newton_backslash_norm=l2.csv};
                \addplot +[very thick] table [x expr={\thisrowno{0}}, y expr={\thisrowno{1}},
                            skip first n=1, col sep=comma]
                {numerical_results/Oscillator/EntropyConservative/Results_HBPC_serial_HB-I2DRK6-3s_relaxation=1_kmax=4_Predictor=Taylor_oscillator_Tend=10_nonlin=1_damped_newton_backslash_norm=l2.csv};
                \addplot +[very thick] table [x expr={\thisrowno{0}}, y expr={\thisrowno{1}},
                            skip first n=1, col sep=comma]
                {numerical_results/Oscillator/EntropyConservative/Results_HBPC_serial_HB-I2DRK6-3s_relaxation=1_kmax=6_Predictor=Taylor_oscillator_Tend=10_nonlin=1_damped_newton_backslash_norm=l2.csv};\addplot +[very thick] table [x expr={\thisrowno{0}}, y expr={\thisrowno{1}},
                            skip first n=1, col sep=comma]
                {numerical_results/Oscillator/EntropyConservative/Results_HBPC_serial_HB-I2DRK6-3s_relaxation=1_kmax=10_Predictor=Taylor_oscillator_Tend=10_nonlin=1_damped_newton_backslash_norm=l2.csv};
            %\legend{$k_{\max}=0$,$k_{\max}=1$,$k_{\max}=2$,$k_{\max}=3$,$k_{\max}=4$,$k_{\max}=6$,$k_{\max}=10$}
            \end{loglogaxis}
            \end{tikzpicture}
        \begin{tikzpicture}[scale=0.6]
            \begin{loglogaxis}[
              cycle list name=hierarchy6,
              xlabel={$\Delta t$},
              ylabel={},
              grid=major,
              yticklabels=\empty,
              legend style={at={(1.1,0.02)},anchor= south east,font=\large},
              title={HBPC(2,8,$k_{\max}$) with relaxation},
              label style={font=\large},
              title style={font=\Large},
              ymax=1e1,ymin=1e-13,
              xticklabel style = {font=\large},
              yticklabel style = {font=\large},
            ]
            \addplot +[very thick] table [x expr={\thisrowno{0}}, y expr={\thisrowno{1}},
                            skip first n=1, col sep=comma]
                {numerical_results/Oscillator/EntropyConservative/Results_HBPC_serial_HB-I2DRK8-4s_relaxation=1_kmax=0_Predictor=Taylor_oscillator_Tend=10_nonlin=1_damped_newton_backslash_norm=l2.csv};
            \addplot +[very thick] table [x expr={\thisrowno{0}}, y expr={\thisrowno{1}},
                            skip first n=1, col sep=comma]
                {numerical_results/Oscillator/EntropyConservative/Results_HBPC_serial_HB-I2DRK8-4s_relaxation=1_kmax=1_Predictor=Taylor_oscillator_Tend=10_nonlin=1_damped_newton_backslash_norm=l2.csv};
                \addplot +[very thick] table [x expr={\thisrowno{0}}, y expr={\thisrowno{1}},
                            skip first n=1, col sep=comma]
                {numerical_results/Oscillator/EntropyConservative/Results_HBPC_serial_HB-I2DRK8-4s_relaxation=1_kmax=2_Predictor=Taylor_oscillator_Tend=10_nonlin=1_damped_newton_backslash_norm=l2.csv};
                \addplot +[very thick] table [x expr={\thisrowno{0}}, y expr={\thisrowno{1}},
                            skip first n=1, col sep=comma]
                {numerical_results/Oscillator/EntropyConservative/Results_HBPC_serial_HB-I2DRK8-4s_relaxation=1_kmax=3_Predictor=Taylor_oscillator_Tend=10_nonlin=1_damped_newton_backslash_norm=l2.csv};
                \addplot +[very thick] table [x expr={\thisrowno{0}}, y expr={\thisrowno{1}},
                            skip first n=1, col sep=comma]
                {numerical_results/Oscillator/EntropyConservative/Results_HBPC_serial_HB-I2DRK8-4s_relaxation=1_kmax=4_Predictor=Taylor_oscillator_Tend=10_nonlin=1_damped_newton_backslash_norm=l2.csv};
                \addplot +[very thick] table [x expr={\thisrowno{0}}, y expr={\thisrowno{1}},
                            skip first n=1, col sep=comma]
                {numerical_results/Oscillator/EntropyConservative/Results_HBPC_serial_HB-I2DRK8-4s_relaxation=1_kmax=6_Predictor=Taylor_oscillator_Tend=10_nonlin=1_damped_newton_backslash_norm=l2.csv};
                \addplot +[very thick] table [x expr={\thisrowno{0}}, y expr={\thisrowno{1}},
                            skip first n=1, col sep=comma]
                {numerical_results/Oscillator/EntropyConservative/Results_HBPC_serial_HB-I2DRK8-4s_relaxation=1_kmax=10_Predictor=Taylor_oscillator_Tend=10_nonlin=1_damped_newton_backslash_norm=l2.csv};
            %\legend{$k_{\max}=0$,$k_{\max}=1$,$k_{\max}=2$,$k_{\max}=3$,$k_{\max}=4$,$k_{\max}=6$,$k_{\max}=10$}
            \end{loglogaxis}
            \end{tikzpicture}
        \begin{tikzpicture}[scale=0.6]
            \begin{loglogaxis}[
              cycle list name=hierarchy6,
              xlabel={$\Delta t$},
              ylabel={},
              grid=major,
              yticklabels=\empty,
              legend style={at={(1.02,0.02)},anchor= south west,font=\normalsize},
              title={HBPC(3,6,$k_{\max}$) with relaxation},
              label style={font=\large},
              title style={font=\Large},
              ymax=1e1,ymin=1e-13,
              xticklabel style = {font=\large},
              yticklabel style = {font=\large}
            ]
            \addplot +[very thick] table [x expr={\thisrowno{0}}, y expr={\thisrowno{1}},
                            skip first n=1, col sep=comma]
                {numerical_results/Oscillator/EntropyConservative/Results_HBPC_serial_HB-I3DRK6-2s_relaxation=1_kmax=0_Predictor=Taylor_oscillator_Tend=10_nonlin=1_damped_newton_backslash_norm=l2.csv};
            \addplot +[very thick] table [x expr={\thisrowno{0}}, y expr={\thisrowno{1}},
                            skip first n=1, col sep=comma]
                {numerical_results/Oscillator/EntropyConservative/Results_HBPC_serial_HB-I3DRK6-2s_relaxation=1_kmax=1_Predictor=Taylor_oscillator_Tend=10_nonlin=1_damped_newton_backslash_norm=l2.csv};
                \addplot +[very thick] table [x expr={\thisrowno{0}}, y expr={\thisrowno{1}},
                            skip first n=1, col sep=comma]
                {numerical_results/Oscillator/EntropyConservative/Results_HBPC_serial_HB-I3DRK6-2s_relaxation=1_kmax=2_Predictor=Taylor_oscillator_Tend=10_nonlin=1_damped_newton_backslash_norm=l2.csv};
                \addplot +[very thick] table [x expr={\thisrowno{0}}, y expr={\thisrowno{1}},
                            skip first n=1, col sep=comma]
                {numerical_results/Oscillator/EntropyConservative/Results_HBPC_serial_HB-I3DRK6-2s_relaxation=1_kmax=3_Predictor=Taylor_oscillator_Tend=10_nonlin=1_damped_newton_backslash_norm=l2.csv};
                \addplot +[very thick] table [x expr={\thisrowno{0}}, y expr={\thisrowno{1}},
                            skip first n=1, col sep=comma]
                {numerical_results/Oscillator/EntropyConservative/Results_HBPC_serial_HB-I3DRK6-2s_relaxation=1_kmax=4_Predictor=Taylor_oscillator_Tend=10_nonlin=1_damped_newton_backslash_norm=l2.csv};
                \addplot +[very thick] table [x expr={\thisrowno{0}}, y expr={\thisrowno{1}},
                            skip first n=1, col sep=comma]
                {numerical_results/Oscillator/EntropyConservative/Results_HBPC_serial_HB-I3DRK6-2s_relaxation=1_kmax=6_Predictor=Taylor_oscillator_Tend=10_nonlin=1_damped_newton_backslash_norm=l2.csv};\addplot +[very thick] table [x expr={\thisrowno{0}}, y expr={\thisrowno{1}},
                            skip first n=1, col sep=comma]
                {numerical_results/Oscillator/EntropyConservative/Results_HBPC_serial_HB-I3DRK6-2s_relaxation=1_kmax=10_Predictor=Taylor_oscillator_Tend=10_nonlin=1_damped_newton_backslash_norm=l2.csv};
            \legend{$k_{\max}=0$,$k_{\max}=1$,$k_{\max}=2$,$k_{\max}=3$,$k_{\max}=4$,$k_{\max}=6$,$k_{\max}=10$}
            \end{loglogaxis}
            \end{tikzpicture}
            %\end{minipage}
        \caption{
        The (serial) HBPC(2, 6, $k_{\max}$) (left), HBPC(2, 8, $k_{\max}$) (middle) and HBPC(3, 6, $\km$) (right)  predictor-corrector scheme of \cite{ZSS22}, see also \cite{SchuetzSealZeifang21}, applied to the entropy-conserving nonlinear oscillator at $\Te = 10$ for various values of $k_{\max}$. Top: without a relaxation procedure. Bottom: with relaxation procedure. The order of convergence to be expected, see Rem.~\ref{rem:order}, is $\min\{6, \km+2\}$ for the HBPC(2, 6, $\km$) scheme, $\min\{8, \km+2\}$ for the HBPC(2, 8, $\km$) scheme and $\min\{6,\km+3\}$ for the HBPC(3, 6, $\km$) scheme. This expected order is met for the unrelaxed version. The relaxed version shows an odd-even decoupling, so for odd $\km$, the order is increased by one.
        }\label{fig:ConvergenceOscillator}
    \end{figure}

\paragraph{Convergence properties}
In a subsequent step, we analyze the convergence properties of the method. Fig.~\ref{fig:ConvergenceOscillator} shows convergence results for two two-derivative and one three-derivative scheme, each with and without relaxation.
From Rem.~\ref{rem:order}, the order of convergence is supposed to be the minimum of $\km + m$ and the ultimate order $q$ of the background Runge-Kutta scheme. It can be seen for the unrelaxed case, that this order is indeed met. For the relaxed version, we see an odd-even decoupling of the order, i.e., for an odd value of $\km$, the order is one order better than expected. The maximum order of consistency, however, remains $q$.
This has been proved in \cite{ranocha2020relaxationHamiltonian} for general B-series methods and the special situation of Euclidean Hamiltonian problems as in this case.
In any case, the error constants seem to be tremendously lower for the relaxed version which is obviously also backed up through the findings from Fig.~\ref{fig:ErrorGrowthOscillator}.

\subsection{Kepler's problem}

To confirm some of the results from the previous section, and to show that the odd-even decoupling is not so much a feature of the method, but more of the underlying problem, we consider here Kepler's problem as in \cite{ranocha2020relaxationHamiltonian}. The problem is given by
\begin{align*}
 \Phi(w) := \begin{pmatrix}
             w_3 \\ w_4 \\ -\frac{w_1}{(w_1^2+w_2^2)^{\frac 3 2}} \\ -\frac{w_2}{(w_1^2+w_2^2)^{\frac 3 2}}
            \end{pmatrix}, \qquad
w(0) = \begin{pmatrix}
        1/2 \\ 0 \\ 0 \\ \sqrt{1/3 }.
       \end{pmatrix}.
\end{align*}
The angular momentum
\begin{align*}
 \eta(w) := w_1w_4 - w_2 w_3
\end{align*}
is a conserved quantity.
For this example, $\dt = 0.5$ is way too coarse, and the relaxed version was not able to run due to the fact that at some point, the relaxation parameter $\gamma$ from \eqref{eq:relaxation} could not be computed anymore. In this way, the relaxed algorithm also gives some extra information on the quality of the solution. Hence, we use smaller $\dt$ here. As in the example before, we start with error growth as a function of $t$ for two values of $\dt$, in this case $\dt = 0.2$ and $\dt = 0.05$, see Fig.~\ref{fig:ErrorGrowthKepler}.
Again, we can see that the error growth for the relaxed method is slower than for the unrelaxed version. It is not a clear linear / quadratic relation as before due to periodic effects, but the overall growth seems in fact to be dominated by linear (relaxed) and quadratic (unrelaxed) terms. Fig.~\ref{fig:ErrorGrowthKepler}, bottom, shows the deviation of the functional $\eta$ from the value $\eta_0 := \eta(w(0))$. As expected, for the relaxed version, it is preserved, while deviations for the unrelaxed algorithm are visible.

\begin{figure} [ht]
        \begin{tikzpicture}[scale=0.6]
            \begin{loglogaxis}[
              cycle list name=hierarchy6,
              xlabel={$\Delta t$},
              ylabel={},
              grid=major,
              legend style={at={(0.02,0.98)},anchor= north west,font=\large},
              title={HBPC(2,6,$k_{\max}$), no relaxation},
              label style={font=\large},
              title style={font=\Large},
              ymax=1e1,ymin=1e-13,
              xticklabel style = {font=\large},
              yticklabel style = {font=\large},
            ]
            \addplot +[very thick] table [x expr={\thisrowno{0}}, y expr={\thisrowno{1}},
                            skip first n=1, col sep=comma]
                {numerical_results/Kepler/EntropyConservative/Results_HBPC_serial_HB-I2DRK6-3s_relaxation=0_kmax=0_Predictor=Taylor_Kepler_Tend=5_damped_newton_backslash_norm=l2.csv};
            \addplot +[very thick] table [x expr={\thisrowno{0}}, y expr={\thisrowno{1}},
                            skip first n=1, col sep=comma]
                {numerical_results/Kepler/EntropyConservative/Results_HBPC_serial_HB-I2DRK6-3s_relaxation=0_kmax=1_Predictor=Taylor_Kepler_Tend=5_damped_newton_backslash_norm=l2.csv};
                \addplot +[very thick] table [x expr={\thisrowno{0}}, y expr={\thisrowno{1}},
                            skip first n=1, col sep=comma]
                {numerical_results/Kepler/EntropyConservative/Results_HBPC_serial_HB-I2DRK6-3s_relaxation=0_kmax=2_Predictor=Taylor_Kepler_Tend=5_damped_newton_backslash_norm=l2.csv};
                \addplot +[very thick] table [x expr={\thisrowno{0}}, y expr={\thisrowno{1}},
                            skip first n=1, col sep=comma]
                {numerical_results/Kepler/EntropyConservative/Results_HBPC_serial_HB-I2DRK6-3s_relaxation=0_kmax=3_Predictor=Taylor_Kepler_Tend=5_damped_newton_backslash_norm=l2.csv};
                \addplot +[very thick] table [x expr={\thisrowno{0}}, y expr={\thisrowno{1}},
                            skip first n=1, col sep=comma]
                {numerical_results/Kepler/EntropyConservative/Results_HBPC_serial_HB-I2DRK6-3s_relaxation=0_kmax=4_Predictor=Taylor_Kepler_Tend=5_damped_newton_backslash_norm=l2.csv};
                \addplot +[very thick] table [x expr={\thisrowno{0}}, y expr={\thisrowno{1}},
                            skip first n=1, col sep=comma]
                {numerical_results/Kepler/EntropyConservative/Results_HBPC_serial_HB-I2DRK6-3s_relaxation=0_kmax=6_Predictor=Taylor_Kepler_Tend=5_damped_newton_backslash_norm=l2.csv};\addplot +[very thick] table [x expr={\thisrowno{0}}, y expr={\thisrowno{1}},
                            skip first n=1, col sep=comma]
                {numerical_results/Kepler/EntropyConservative/Results_HBPC_serial_HB-I2DRK6-3s_relaxation=0_kmax=10_Predictor=Taylor_Kepler_Tend=5_damped_newton_backslash_norm=l2.csv};
            %\legend{$k_{\max}=0$,$k_{\max}=1$,$k_{\max}=2$,$k_{\max}=3$,$k_{\max}=4$,$k_{\max}=6$,$k_{\max}=10$}
            \end{loglogaxis}
            \end{tikzpicture}
            \begin{tikzpicture}[scale=0.6]
            \begin{loglogaxis}[
              cycle list name=hierarchy6,
              xlabel={$\Delta t$},
              ylabel={},
              grid=major,
              legend style={at={(0.02,0.98)},anchor= north west,font=\large},
              title={HBPC(2,8,$k_{\max})$, no relaxation},
              label style={font=\large},
              title style={font=\Large},
              yticklabels=\empty,
              ymax=1e1,ymin=1e-13,
              xticklabel style = {font=\large},
              yticklabel style = {font=\large},
            ]
            \addplot +[very thick] table [x expr={\thisrowno{0}}, y expr={\thisrowno{1}},
                            skip first n=1, col sep=comma]
                {numerical_results/Kepler/EntropyConservative/Results_HBPC_serial_HB-I2DRK8-4s_relaxation=0_kmax=0_Predictor=Taylor_Kepler_Tend=5_damped_newton_backslash_norm=l2.csv};
            \addplot +[very thick] table [x expr={\thisrowno{0}}, y expr={\thisrowno{1}},
                            skip first n=1, col sep=comma]
                {numerical_results/Kepler/EntropyConservative/Results_HBPC_serial_HB-I2DRK8-4s_relaxation=0_kmax=1_Predictor=Taylor_Kepler_Tend=5_damped_newton_backslash_norm=l2.csv};
                \addplot +[very thick] table [x expr={\thisrowno{0}}, y expr={\thisrowno{1}},
                            skip first n=1, col sep=comma]
                {numerical_results/Kepler/EntropyConservative/Results_HBPC_serial_HB-I2DRK8-4s_relaxation=0_kmax=2_Predictor=Taylor_Kepler_Tend=5_damped_newton_backslash_norm=l2.csv};
                \addplot +[very thick] table [x expr={\thisrowno{0}}, y expr={\thisrowno{1}},
                            skip first n=1, col sep=comma]
                {numerical_results/Kepler/EntropyConservative/Results_HBPC_serial_HB-I2DRK8-4s_relaxation=0_kmax=3_Predictor=Taylor_Kepler_Tend=5_damped_newton_backslash_norm=l2.csv};
                \addplot +[very thick] table [x expr={\thisrowno{0}}, y expr={\thisrowno{1}},
                            skip first n=1, col sep=comma]
                {numerical_results/Kepler/EntropyConservative/Results_HBPC_serial_HB-I2DRK8-4s_relaxation=0_kmax=4_Predictor=Taylor_Kepler_Tend=5_damped_newton_backslash_norm=l2.csv};
                \addplot +[very thick] table [x expr={\thisrowno{0}}, y expr={\thisrowno{1}},
                            skip first n=1, col sep=comma]
                {numerical_results/Kepler/EntropyConservative/Results_HBPC_serial_HB-I2DRK8-4s_relaxation=0_kmax=6_Predictor=Taylor_Kepler_Tend=5_damped_newton_backslash_norm=l2.csv};\addplot +[very thick] table [x expr={\thisrowno{0}}, y expr={\thisrowno{1}},
                            skip first n=1, col sep=comma]
                {numerical_results/Kepler/EntropyConservative/Results_HBPC_serial_HB-I2DRK8-4s_relaxation=0_kmax=10_Predictor=Taylor_Kepler_Tend=5_damped_newton_backslash_norm=l2.csv};
            %\legend{$k_{\max}=0$,$k_{\max}=1$,$k_{\max}=2$,$k_{\max}=3$,$k_{\max}=4$,$k_{\max}=6$,$k_{\max}=10$}
            \end{loglogaxis}
            \end{tikzpicture}
            \begin{tikzpicture}[scale=0.6]
            \begin{loglogaxis}[
              cycle list name=hierarchy6,
              xlabel={$\Delta t$},
              ylabel={},
              grid=major,
              legend style={at={(1.02,0.02)},anchor= south west,font=\normalsize},
              title={HBPC(3,6,$k_{\max})$, no relaxation},
              label style={font=\large},
              title style={font=\Large},
              yticklabels=\empty,
              ymax=1e1,ymin=1e-13,
              xticklabel style = {font=\large},
              yticklabel style = {font=\large},
            ]
            \addplot +[very thick] table [x expr={\thisrowno{0}}, y expr={\thisrowno{1}},
                            skip first n=1, col sep=comma]
                {numerical_results/Kepler/EntropyConservative/Results_HBPC_serial_HB-I3DRK6-2s_relaxation=0_kmax=0_Predictor=Taylor_Kepler_Tend=5_damped_newton_backslash_norm=l2.csv};
            \addplot +[very thick] table [x expr={\thisrowno{0}}, y expr={\thisrowno{1}},
                            skip first n=1, col sep=comma]
                {numerical_results/Kepler/EntropyConservative/Results_HBPC_serial_HB-I3DRK6-2s_relaxation=0_kmax=1_Predictor=Taylor_Kepler_Tend=5_damped_newton_backslash_norm=l2.csv};
                \addplot +[very thick] table [x expr={\thisrowno{0}}, y expr={\thisrowno{1}},
                            skip first n=1, col sep=comma]
                {numerical_results/Kepler/EntropyConservative/Results_HBPC_serial_HB-I3DRK6-2s_relaxation=0_kmax=2_Predictor=Taylor_Kepler_Tend=5_damped_newton_backslash_norm=l2.csv};
                \addplot +[very thick] table [x expr={\thisrowno{0}}, y expr={\thisrowno{1}},
                            skip first n=1, col sep=comma]
                {numerical_results/Kepler/EntropyConservative/Results_HBPC_serial_HB-I3DRK6-2s_relaxation=0_kmax=3_Predictor=Taylor_Kepler_Tend=5_damped_newton_backslash_norm=l2.csv};
                \addplot +[very thick] table [x expr={\thisrowno{0}}, y expr={\thisrowno{1}},
                            skip first n=1, col sep=comma]
                {numerical_results/Kepler/EntropyConservative/Results_HBPC_serial_HB-I3DRK6-2s_relaxation=0_kmax=4_Predictor=Taylor_Kepler_Tend=5_damped_newton_backslash_norm=l2.csv};
                \addplot +[very thick] table [x expr={\thisrowno{0}}, y expr={\thisrowno{1}},
                            skip first n=1, col sep=comma]
                {numerical_results/Kepler/EntropyConservative/Results_HBPC_serial_HB-I3DRK6-2s_relaxation=0_kmax=6_Predictor=Taylor_Kepler_Tend=5_damped_newton_backslash_norm=l2.csv};
                \addplot +[very thick] table [x expr={\thisrowno{0}}, y expr={\thisrowno{1}},
                            skip first n=1, col sep=comma]
                {numerical_results/Kepler/EntropyConservative/Results_HBPC_serial_HB-I3DRK6-2s_relaxation=0_kmax=10_Predictor=Taylor_Kepler_Tend=5_damped_newton_backslash_norm=l2.csv};
            \legend{$k_{\max}=0$,$k_{\max}=1$,$k_{\max}=2$,$k_{\max}=3$,$k_{\max}=4$,$k_{\max}=6$,$k_{\max}=10$}
            \end{loglogaxis}
            \end{tikzpicture} \\
        \begin{tikzpicture}[scale=0.6]
            \begin{loglogaxis}[
              cycle list name=hierarchy6,
              xlabel={$\Delta t$},
              ylabel={},
              grid=major,
              title={HBPC(2,6,$k_{\max}$) with relaxation},
              label style={font=\large},
              title style={font=\Large},
              ymax=1e1,ymin=1e-13,
              xticklabel style = {font=\large},
              yticklabel style = {font=\large},
            ]
            \addplot +[very thick] table [x expr={\thisrowno{0}}, y expr={\thisrowno{1}},
                            skip first n=1, col sep=comma]
                {numerical_results/Kepler/EntropyConservative/Results_HBPC_serial_HB-I2DRK6-3s_relaxation=1_kmax=0_Predictor=Taylor_Kepler_Tend=5_damped_newton_backslash_norm=l2.csv};
            \addplot +[very thick] table [x expr={\thisrowno{0}}, y expr={\thisrowno{1}},
                            skip first n=1, col sep=comma]
                {numerical_results/Kepler/EntropyConservative/Results_HBPC_serial_HB-I2DRK6-3s_relaxation=1_kmax=1_Predictor=Taylor_Kepler_Tend=5_damped_newton_backslash_norm=l2.csv};
                \addplot +[very thick] table [x expr={\thisrowno{0}}, y expr={\thisrowno{1}},
                            skip first n=1, col sep=comma]
                {numerical_results/Kepler/EntropyConservative/Results_HBPC_serial_HB-I2DRK6-3s_relaxation=1_kmax=2_Predictor=Taylor_Kepler_Tend=5_damped_newton_backslash_norm=l2.csv};
                \addplot +[very thick] table [x expr={\thisrowno{0}}, y expr={\thisrowno{1}},
                            skip first n=1, col sep=comma]
                {numerical_results/Kepler/EntropyConservative/Results_HBPC_serial_HB-I2DRK6-3s_relaxation=1_kmax=3_Predictor=Taylor_Kepler_Tend=5_damped_newton_backslash_norm=l2.csv};
                \addplot +[very thick] table [x expr={\thisrowno{0}}, y expr={\thisrowno{1}},
                            skip first n=1, col sep=comma]
                {numerical_results/Kepler/EntropyConservative/Results_HBPC_serial_HB-I2DRK6-3s_relaxation=1_kmax=4_Predictor=Taylor_Kepler_Tend=5_damped_newton_backslash_norm=l2.csv};
                \addplot +[very thick] table [x expr={\thisrowno{0}}, y expr={\thisrowno{1}},
                            skip first n=1, col sep=comma]
                {numerical_results/Kepler/EntropyConservative/Results_HBPC_serial_HB-I2DRK6-3s_relaxation=1_kmax=6_Predictor=Taylor_Kepler_Tend=5_damped_newton_backslash_norm=l2.csv};\addplot +[very thick] table [x expr={\thisrowno{0}}, y expr={\thisrowno{1}},
                            skip first n=1, col sep=comma]
                {numerical_results/Kepler/EntropyConservative/Results_HBPC_serial_HB-I2DRK6-3s_relaxation=1_kmax=10_Predictor=Taylor_Kepler_Tend=5_damped_newton_backslash_norm=l2.csv};
            %\legend{$k_{\max}=0$,$k_{\max}=1$,$k_{\max}=2$,$k_{\max}=3$,$k_{\max}=4$,$k_{\max}=6$,$k_{\max}=10$}
            \end{loglogaxis}
            \end{tikzpicture}
        \begin{tikzpicture}[scale=0.6]
            \begin{loglogaxis}[
              cycle list name=hierarchy6,
              xlabel={$\Delta t$},
              ylabel={},
              grid=major,
              yticklabels=\empty,
              title={HBPC(2,8,$k_{\max}$) with relaxation},
              label style={font=\large},
              title style={font=\Large},
              ymax=1e1,ymin=1e-13,
              xticklabel style = {font=\large},
              yticklabel style = {font=\large},
            ]
            \addplot +[very thick] table [x expr={\thisrowno{0}}, y expr={\thisrowno{1}},
                            skip first n=1, col sep=comma]
                {numerical_results/Kepler/EntropyConservative/Results_HBPC_serial_HB-I2DRK8-4s_relaxation=1_kmax=0_Predictor=Taylor_Kepler_Tend=5_damped_newton_backslash_norm=l2.csv};
            \addplot +[very thick] table [x expr={\thisrowno{0}}, y expr={\thisrowno{1}},
                            skip first n=1, col sep=comma]
                {numerical_results/Kepler/EntropyConservative/Results_HBPC_serial_HB-I2DRK8-4s_relaxation=1_kmax=1_Predictor=Taylor_Kepler_Tend=5_damped_newton_backslash_norm=l2.csv};
                \addplot +[very thick] table [x expr={\thisrowno{0}}, y expr={\thisrowno{1}},
                            skip first n=1, col sep=comma]
                {numerical_results/Kepler/EntropyConservative/Results_HBPC_serial_HB-I2DRK8-4s_relaxation=1_kmax=2_Predictor=Taylor_Kepler_Tend=5_damped_newton_backslash_norm=l2.csv};
                \addplot +[very thick] table [x expr={\thisrowno{0}}, y expr={\thisrowno{1}},
                            skip first n=1, col sep=comma]
                {numerical_results/Kepler/EntropyConservative/Results_HBPC_serial_HB-I2DRK8-4s_relaxation=1_kmax=3_Predictor=Taylor_Kepler_Tend=5_damped_newton_backslash_norm=l2.csv};
                \addplot +[very thick] table [x expr={\thisrowno{0}}, y expr={\thisrowno{1}},
                            skip first n=1, col sep=comma]
                {numerical_results/Kepler/EntropyConservative/Results_HBPC_serial_HB-I2DRK8-4s_relaxation=1_kmax=4_Predictor=Taylor_Kepler_Tend=5_damped_newton_backslash_norm=l2.csv};
                \addplot +[very thick] table [x expr={\thisrowno{0}}, y expr={\thisrowno{1}},
                            skip first n=1, col sep=comma]
                {numerical_results/Kepler/EntropyConservative/Results_HBPC_serial_HB-I2DRK8-4s_relaxation=1_kmax=6_Predictor=Taylor_Kepler_Tend=5_damped_newton_backslash_norm=l2.csv};
                \addplot +[very thick] table [x expr={\thisrowno{0}}, y expr={\thisrowno{1}},
                            skip first n=1, col sep=comma]
                {numerical_results/Kepler/EntropyConservative/Results_HBPC_serial_HB-I2DRK8-4s_relaxation=1_kmax=10_Predictor=Taylor_Kepler_Tend=5_damped_newton_backslash_norm=l2.csv};
            %\legend{$k_{\max}=0$,$k_{\max}=1$,$k_{\max}=2$,$k_{\max}=3$,$k_{\max}=4$,$k_{\max}=6$,$k_{\max}=10$}
            \end{loglogaxis}
            \end{tikzpicture}
        \begin{tikzpicture}[scale=0.6]
            \begin{loglogaxis}[
              cycle list name=hierarchy6,
              xlabel={$\Delta t$},
              ylabel={},
              grid=major,
              yticklabels=\empty,
              legend style={at={(1.02,0.02)},anchor= south west,font=\normalsize},
              title={HBPC(3,6,$k_{\max}$) with relaxation},
              label style={font=\large},
              title style={font=\Large},
              ymax=1e1,ymin=1e-13,
              xticklabel style = {font=\large},
              yticklabel style = {font=\large},
            ]
            \addplot +[very thick] table [x expr={\thisrowno{0}}, y expr={\thisrowno{1}},
                            skip first n=1, col sep=comma]
                {numerical_results/Kepler/EntropyConservative/Results_HBPC_serial_HB-I3DRK6-2s_relaxation=1_kmax=0_Predictor=Taylor_Kepler_Tend=5_damped_newton_backslash_norm=l2.csv};
            \addplot +[very thick] table [x expr={\thisrowno{0}}, y expr={\thisrowno{1}},
                            skip first n=1, col sep=comma]
                {numerical_results/Kepler/EntropyConservative/Results_HBPC_serial_HB-I3DRK6-2s_relaxation=1_kmax=1_Predictor=Taylor_Kepler_Tend=5_damped_newton_backslash_norm=l2.csv};
                \addplot +[very thick] table [x expr={\thisrowno{0}}, y expr={\thisrowno{1}},
                            skip first n=1, col sep=comma]
                {numerical_results/Kepler/EntropyConservative/Results_HBPC_serial_HB-I3DRK6-2s_relaxation=1_kmax=2_Predictor=Taylor_Kepler_Tend=5_damped_newton_backslash_norm=l2.csv};
                \addplot +[very thick] table [x expr={\thisrowno{0}}, y expr={\thisrowno{1}},
                            skip first n=1, col sep=comma]
                {numerical_results/Kepler/EntropyConservative/Results_HBPC_serial_HB-I3DRK6-2s_relaxation=1_kmax=3_Predictor=Taylor_Kepler_Tend=5_damped_newton_backslash_norm=l2.csv};
                \addplot +[very thick] table [x expr={\thisrowno{0}}, y expr={\thisrowno{1}},
                            skip first n=1, col sep=comma]
                {numerical_results/Kepler/EntropyConservative/Results_HBPC_serial_HB-I3DRK6-2s_relaxation=1_kmax=4_Predictor=Taylor_Kepler_Tend=5_damped_newton_backslash_norm=l2.csv};
                \addplot +[very thick] table [x expr={\thisrowno{0}}, y expr={\thisrowno{1}},
                            skip first n=1, col sep=comma]
                {numerical_results/Kepler/EntropyConservative/Results_HBPC_serial_HB-I3DRK6-2s_relaxation=1_kmax=6_Predictor=Taylor_Kepler_Tend=5_damped_newton_backslash_norm=l2.csv};\addplot +[very thick] table [x expr={\thisrowno{0}}, y expr={\thisrowno{1}},
                            skip first n=1, col sep=comma]
                {numerical_results/Kepler/EntropyConservative/Results_HBPC_serial_HB-I3DRK6-2s_relaxation=1_kmax=10_Predictor=Taylor_Kepler_Tend=5_damped_newton_backslash_norm=l2.csv};
            \legend{$k_{\max}=0$,$k_{\max}=1$,$k_{\max}=2$,$k_{\max}=3$,$k_{\max}=4$,$k_{\max}=6$,$k_{\max}=10$}
            \end{loglogaxis}
            \end{tikzpicture}
        \caption{The (serial) HBPC(2, 6, $k_{\max}$) (left), HBPC(2, 8, $k_{\max}$) (middle) and HBPC(3, 6, $\km$) (right)  predictor-corrector scheme of \cite{ZSS22}, see also \cite{SchuetzSealZeifang21}, applied to Kepler's problem at $\Te = 5$ for various values of $k_{\max}$. Top: without a relaxation procedure. Bottom: with relaxation procedure. The order of convergence to be expected, see Rem.~\ref{rem:order}, is $\min\{6, \km+2\}$ for the HBPC(2, 6, $\km$) scheme, $\min\{8, \km+2\}$ for the HBPC(2, 8, $\km$) scheme and $\min\{6,\km+3\}$ for the HBPC(3, 6, $\km$) scheme. This expected order is met for both the relaxed and the unrelaxed version. In contrast to the oscillator problem, see Fig.~\ref{fig:ConvergenceOscillator}, only for $\km = 1$ and $\km = 2$, one can see significant differences in the error.
        Please note that the reference solution against which we compute the numerical error is also computed numerically (with another scheme) on a fine resolution. This explains why at about an error level of $10^{-11}$, convergence stalls for all combinations.
        }\label{fig:ConvergenceKepler}
    \end{figure}
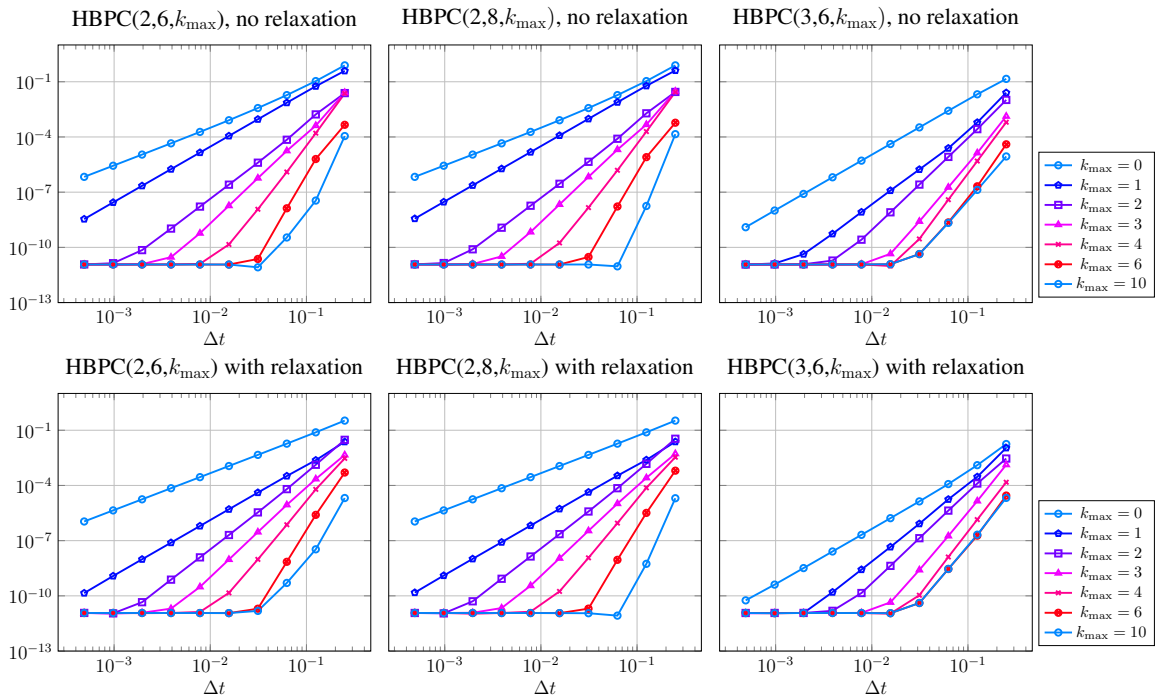

Fig.~\ref{fig:ConvergenceKepler} shows convergence plots for the three different methods used here, two two-derivative and one three-derivative method. In contrast to the results before, there is no odd-even decoupling anymore, and the order of convergence of $\min\{\km+m,q\}$ is clearly met. This clearly indicates that this odd-even decoupling of the order for the relaxed version cannot be expected for all testcases, and is really a feature of the previous problem. Also the reduction of the error constant is only visible for $\km = 1$ (here it is the most prominent) and for $\km = 2$ (slightly). For the higher $\km$, this effect is not really significant.

	\begin{figure} [ht]
        \centering
        \sidecaption
        \begin{minipage}{0.6\textwidth}
        \begin{tikzpicture}[scale=0.6]
            \begin{axis}[
              xlabel={Time $t$},
              ylabel={Numerical error},
              grid=major,
              legend style={at={(0.02,0.98)},anchor= north west,font=\large},
              title={Error growth for $\dt = 0.2$},
              label style={font=\large},
              title style={font=\Large, xshift=1.5ex},
              xticklabel style = {font=\large},
              yticklabel style = {font=\large},
            ]
            \addplot +[very thick,mark=none,dashed] table [x expr={\thisrowno{0}}, y expr={\thisrowno{1}}]
                {numerical_results/Kepler/EntropyConservative/ErrorGrowth_HBPC_serial_HB-I2DRK6-3s_relaxation=1_kmax=4_Predictor=Taylor_kepler_Tend=10_damped_newton_backslash_Nt=50_norm=l2.csv};
            \addplot +[very thick,mark=none] table [x expr={\thisrowno{0}}, y expr={\thisrowno{1}}]
                {numerical_results/Kepler/EntropyConservative/ErrorGrowth_HBPC_serial_HB-I2DRK6-3s_relaxation=0_kmax=4_Predictor=Taylor_kepler_Tend=10_damped_newton_backslash_Nt=50_norm=l2.csv};
                \legend{{With relaxation}, {Without relaxation}}
            \end{axis}
            \end{tikzpicture}
        \begin{tikzpicture}[scale=0.6]
            \begin{axis}[
              xlabel={Time $t$},
              grid=major,
              legend style={at={(0.02,0.98)},anchor= north west,font=\large},
              title={Error growth for $\dt = 0.05$},
              label style={font=\large},
              title style={font=\Large,xshift=2.5ex},
              xticklabel style = {font=\large},
              yticklabel style = {font=\large},
            ]
            \addplot +[very thick,mark=none,dashed] table [x expr={\thisrowno{0}}, y expr={\thisrowno{1}}]
                {numerical_results/Kepler/EntropyConservative/ErrorGrowth_HBPC_serial_HB-I2DRK6-3s_relaxation=1_kmax=4_Predictor=Taylor_kepler_Tend=10_damped_newton_backslash_Nt=200_norm=l2.csv};
            \addplot +[very thick,mark=none] table [x expr={\thisrowno{0}}, y expr={\thisrowno{1}}]
                {numerical_results/Kepler/EntropyConservative/ErrorGrowth_HBPC_serial_HB-I2DRK6-3s_relaxation=0_kmax=4_Predictor=Taylor_kepler_Tend=10_damped_newton_backslash_Nt=200_norm=l2.csv};
                \legend{{With relaxation}, {Without relaxation}}
            \end{axis}
            \end{tikzpicture} \\
            \begin{tikzpicture}[scale=0.6]
            \begin{axis}[
              xlabel={Time $t$},
              ylabel={$\eta(w^n)-\eta(w(0))$},
              grid=major,
              legend style={at={(0.02,0.02)},anchor= south west,font=\large},
              title={Evolution of $\eta-\eta_0$ for $\dt = 0.2$},
              label style={font=\large},
              y label style={yshift=1.5ex},
              title style={font=\Large},
              xticklabel style = {font=\large},
              yticklabel style={/pgf/number format/fixed, /pgf/number format/precision=6},
            ]
            \addplot +[very thick,mark=none,dashed] table [x expr={\thisrowno{0}}, y expr={\thisrowno{2}-0.28867513}]
                {numerical_results/Kepler/EntropyConservative/ErrorGrowth_HBPC_serial_HB-I2DRK6-3s_relaxation=1_kmax=4_Predictor=Taylor_kepler_Tend=10_damped_newton_backslash_Nt=50_norm=l2.csv};
            \addplot +[very thick,mark=none] table [x expr={\thisrowno{0}}, y expr={\thisrowno{2}-0.28867513}]
                {numerical_results/Kepler/EntropyConservative/ErrorGrowth_HBPC_serial_HB-I2DRK6-3s_relaxation=0_kmax=4_Predictor=Taylor_kepler_Tend=10_damped_newton_backslash_Nt=50_norm=l2.csv};
                \legend{{With relaxation}, {Without relaxation}}
            \end{axis}
            \end{tikzpicture}
        \begin{tikzpicture}[scale=0.6]
            \begin{axis}[
              xlabel={Time $t$},%ylabel={$\eta(w^n)$},
              grid=major,
              legend style={at={(0.98,0.02)},anchor= south east,font=\large},
              title={Evolution of $\eta-\eta_0$ for $\dt = 0.05$},
              label style={font=\large},
              title style={font=\Large},
              xticklabel style = {font=\large},
              yticklabel style={/pgf/number format/fixed, /pgf/number format/precision=10},
            ]
            \addplot +[very thick,mark=none,dashed] table [x expr={\thisrowno{0}}, y expr={\thisrowno{2}-0.28867513}]
                {numerical_results/Kepler/EntropyConservative/ErrorGrowth_HBPC_serial_HB-I2DRK6-3s_relaxation=1_kmax=4_Predictor=Taylor_kepler_Tend=10_damped_newton_backslash_Nt=200_norm=l2.csv};
            \addplot +[very thick,mark=none] table [x expr={\thisrowno{0}}, y expr={\thisrowno{3}}]
                {numerical_results/Kepler/EntropyConservative/ErrorGrowth_HBPC_serial_HB-I2DRK6-3s_relaxation=0_kmax=4_Predictor=Taylor_kepler_Tend=10_damped_newton_backslash_Nt=200_norm=l2.csv};
                \legend{{With relaxation}, {Without relaxation}}
            \end{axis}
            \end{tikzpicture}
            \end{minipage}
        \caption{The (serial) HBPC(2, 6, 4) predictor-corrector scheme of \cite{ZSS22}, see also \cite{SchuetzSealZeifang21}, applied to Kepler's problem with $\Te = 10$. Left are numerical results for $\dt = 0.2$, right are results for $\dt = 0.05$. Top: numerical error as a function of time; bottom: the deviation in the functional $\eta$ evaluated for the discrete solution. It is clearly visible that in all cases, the relaxed method behaves significantly better than its unrelaxed counterpart.
        %Please note: For $\dt = 0.05$ and the unrelaxed version, the error in $\eta$ is $\O(10^{-8})$, which is why it is not visible in the plot anymore.
        }\label{fig:ErrorGrowthKepler}
    \end{figure}
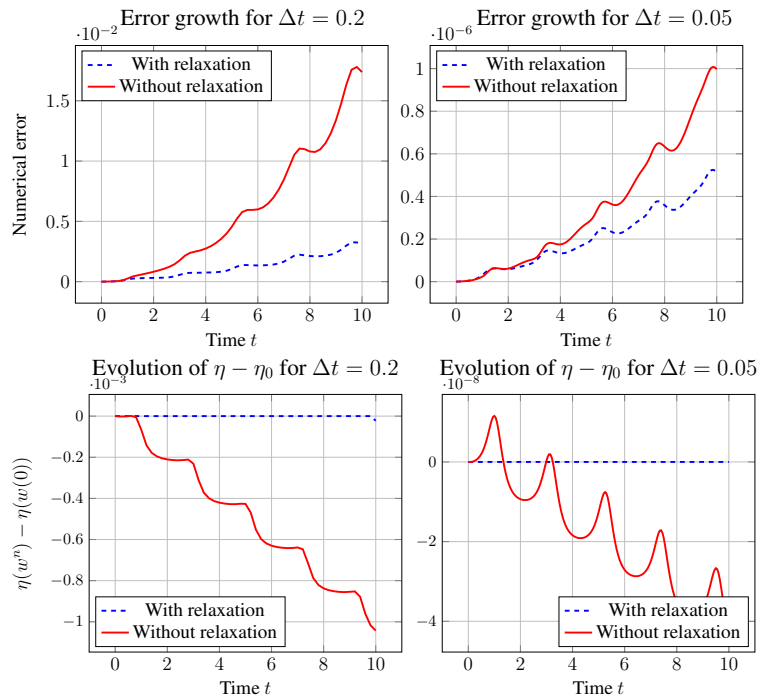

\section{Conclusion and outlook}\label{sec:conout}
In this paper, we have combined recently developed relaxation techniques with also rather recently developed predictor-corrector time integration schemes. It has been shown that this can reduce error constants, and preserve functionals even if the general error level is high.

Obviously, many things are left to do. Currently, we are analyzing, both numerically and analytically, the combination of very general multiderivative methods and relaxation
\begin{itemize}
 \item with respect to convergence properties for many different test problems, including suitably discretized PDEs,
 \item with respect to stability, in particular whether relaxation can change A- and L-stability properties of given methods,
 \item with respect to existence of $\gamma$ and order considerations.
\end{itemize}
Also, dissipative problems, i.e., problems where, contrary to Eq.~\eqref{eq:conservedfunctional}, the functional $\eta$ is not preserved, but decreases over time, i.e., where there holds $\frac{\mathrm{d}}{\mathrm d t} \eta(y(t)) \leq 0$, are subject to investigation.

\begin{acknowledgement}
HR was supported by the Deutsche Forschungsgemeinschaft
(DFG, German Research Foundation, project number 513301895)
and the Daimler und Benz Stiftung (Daimler and Benz foundation,
project number 32-10/22). ET was funded by the Fonds voor Wetenschappelijk Onderzoek (FWO, Belgium) - project no. G052419N.

\end{acknowledgement}

%%%%%%%%%%%%%%%%%%%%%%%%%%%%%%

% Use this code if you wish to generate your bibliography with BibTeX;
% please replace first the string "demo" below with the name(s) of
% the BibTeX data base(s) you want to use.
% !!!!!!!!!!!!!!!!!!!!!!!!!!!!!!!!!!!!!!!!!!!!!!!!!!!!!!!!!!!!!!!!!!!!!!!!!!!!!!!!!
% !!! The resulting bibliography-output (the contents of the .bbl file)
% !!! must be pasted into this file before submission.
% !!!!!!!!!!!!!!!!!!!!!!!!!!!!!!!!!!!!!!!!!!!!!!!!!!!!!!!!!!!!!!!!!!!!!!!!!!!!!!!!!
%\todo{The contents of the .bbl file must be pasted into this file before submission.} %TODO
%\bibliographystyle{pamm}
%\bibliography{ListPaper.bib}

\providecommand{\WileyBibTextsc}{}
\let\textsc\WileyBibTextsc
\providecommand{\othercit}{}
\providecommand{\jr}[1]{#1}
\providecommand{\etal}{~et~al.}

%
% Replace the following example bibliography with your references
% before submission:
% \begin{thebibliography}{1}
%
% \bibitem{bib1}%
%  F.\,M. Firstauthorfamilyname, F.\,M. Secondauthorfamilyname, and
%   C.~Lastauthorfamilyname,
%  Abbreviatedjournalname \textbf{volume}, page (year).
%
% \bibitem{bib2}%
%  F.~Examplename and  I.\,E. Anotherauthorname,
%  phys. stat. sol. (a) \textbf{1}, 111 (2050).
%
% \bibitem{bib3}%
%  A.~Firstauthorname,  B.~Secondauthorname,  and
%   C.~Thirdauthorname,
% Here Goes the Title of the Book (Publisher, City, year), p.\,111.
%
% \bibitem{bib4}%
%  A.~Firsteditorname,  B.~Secondeditorname,  and
%   C.~Thirdeditorname (eds.),
% Here Goes the Title of the Edited Book (Wiley-VCH, Berlin, 2050), p.\,111.
%
% \bibitem{bib5}%
%  D.~Contributorname,
%  in: The Title of the Book, edited by The Name of the Editors, Followed by
%   the Title of the Series of Books (Publisher, City, year), chap.~1.
%
% \bibitem{bib6}%
%  A.~Lastbutnotleastname,
%  Proceedings 1st Dummy Conference on Citation Formatting, City,
%   Country, Part A (Publisher, City, year),  pp.\,1--11.
%
% \end{thebibliography}
\end{document}